\documentclass[10pt]{amsart}
\usepackage{latexsym}
\usepackage{amssymb}
\usepackage{amsmath}
\usepackage{mathrsfs}
\usepackage{graphicx}
\font\eightsc=cmcsc8

\begin{document}

\title{
Excursus into the History of Calculus
}

\author{S.~S. Kutateladze}
\address[]{
Sobolev Institute of Mathematics\newline
\indent 4 Koptyug Avenue\newline
\indent Novosibirsk, 630090\newline
\indent Russia}
\email{
sskut@member.ams.org
}
\begin{abstract}
This is a brief overview of some turning points
in the history of infinitesimals.
\end{abstract}
\thanks{The Russian version of this talk appeared firstly in
the mimeographed notes ``Fundamentals of Nonstandard Mathematical Analysis'' for the students of
Novosibirsk State University in~1984. Its English versions served
as introduction to {\it Nonstandard Methods of Analysis}
by A.~Kusraev and S.~Kutateladze
(Kluwer Academic Publisher, 1994) and {\it Infinitesimal Analysis} by E.~Gordon, A.~Kusraev,
and S.~Kutateladze (Kluwer Academic Publishers, 2002)}

\maketitle

The ideas of differential and integral calculus are traceable
from  the remote ages,  intertwining  tightly with
the most fundamental mathematical concepts.

I admit readily
that
to present   the evolution
of views of mathematical objects and   the history
of the processes of  calculation
and measurement  which gave an impetus to  the modern theory of
infinitesimals  requires  the Herculean efforts  far beyond
my abilities and intentions.

The matter is significantly aggravated by the fact that
the history of mathematics has always fallen  victim to the notorious
incessant attempts at providing an apologia for
all stylish brand-new conceptions and misconceptions.
In particular, many available expositions of the evolution
of calculus
could hardly be praised as complete, fair, and unbiased.
One-sided views of the nature of the differential and the integral,
hypertrophy of the role of the limit
and neglect of the infinitesimal
have been  spread  so widely in the recent decades that
it is impossible to ignore their existence.

It has become a~truism to say (cf. ~\cite{CourRob}):

\begin{itemize}
\item{}{\small
The genuine foundations of analysis
have  for a~long time been surrounded with mystery as a~result
of unwillingness to admit that the notion of limit enjoys
an~exclusive right to be the~source of new methods.}
\end{itemize}

\noindent
However,
Pontryagin was right to remark  in \cite[pp.~64--65]{Pont} that

\begin{itemize}
\item{}{\small
In a~historical sense, integral and differential
calculus had already been among the established areas
of mathematics long before the theory of limits.
The latter originated as ~superstructure over an existent theory.
Many physicists opine that the so-called rigorous definitions
of derivative and integral are in no way necessary
for satisfactory comprehension of differential and integral calculus.
I share this viewpoint.}
\end{itemize}

\noindent
Considering the above,  it is worthwhile to
discuss a few turning points and crucial ideas
in the  evolution of analysis as expressed in the  words
of classics.
The choice of the corresponding fragments is doomed to be subjective.
Nevertheless, the selection below  seems sufficient
for anyone to acquire a~critical attitude to the numerous incomplete and misleading
delineations of the evolution of infinitesimal methods.

\section{G.~W.~Leibniz and  I.~Newton}

The ancient name for differential and integral calculus is
``infinitesimal analysis.''\footnote   {This term
was used in 1748 by Leonhard Euler in
{\it Introductio in Analysin Infinitorum\/}
\cite{Euler} (cf. \cite[p.~324]{Kline}).}

The first textbook on
this subject was published as far back as 1696
under the title  {\it Analyse des infiniment petits pour
l'intelligence des lignes courbe}.
The textbook was compiled by de l'H\^opital
as a~result of his
contacts with J.~Bernoulli
(senior), one of the most famous
disciples of Leibniz.

The history of  creation of mathematical analysis,
the scientific legacy
of its founders and  their
personal relations have been studied in full detail and even
scrutinized.
Each fair attempt is welcome at reconstructing the train of thought
of the men of genius
and
elucidating the ways to new knowledge and keen vision.
We must however bear
in mind the principal differences between draft
papers and notes, personal letters to colleagues, and the articles
written especially for publication.
It is therefore reasonable to look at the
``official''  presentation of  ~Leibniz's and
~Newton's views of infinitesimals.

The first publication on differential calculus  was
Leibniz's article
``Nova methodus pro maximis et minimis,
itemque tangentibus, quae nec fractals nec
irrationales quantitates moratur,
et singulare pro illis calculi genus''
(see \cite{Leib}).
This article was published in the Leipzig journal ``Acta Eruditorum''
more than three centuries ago in 1684.

Leibniz gave the
following definition of differential.
Considering
a curve $YY$ and a tangent at a~ fixed point $Y$ on the curve
which corresponds to a coordinate $X$ on the axis $AX$ and denoting
by $D$ the  intersection point  of the tangent and  axis, Leibniz
wrote:

\begin{itemize}
\item{}{\small
Now some straight line selected arbitrarily is called
$dx$ and another line
whose ratio to $dx$ is the same as  of
$\dots$  $y$ $\dots$  to $XD$
is called $\dots$  $dy$ or difference ({\it  differentia}) $\dots$
of $y$ $\dots$.}
\end{itemize}

\noindent
The essential details of the picture accompanying
this text are reproduced in~ Fig.~1.

By Leibniz,
given an~arbitrary
$dx$ and considering  the function   $x\mapsto y(x)$ at
a~point ~$x$, we obtain
$$
dy:=\frac{YX}{XD}dx.
$$
In other words, the differential
of a function is defined as the appropriate
linear mapping in the manner fully acceptable  to the majority
of the today's teachers of analysis.

Leibniz
was a~deep thinker and polymath who believed (see \cite[pp.~492--493]{Leib2}) that

\begin{itemize}
\item{}{\small
the invention of the syllogistic form ranks among the most beautiful and
even the most important discoveries of the human mind.
This is a~sort of
{\it universal mathematics\/}
whose significance  has not yet been completely comprehended.
It can be said to incarnate the art of faultlessness ...\,.}
\end{itemize}

\noindent
Leibniz understood definitely that the description and substantiation of
the
algorithm of differential calculus
(in that way he
referred to the rules of differentiation)
required clarifying the concept of tangent.
He proceeded with explaining that

\begin{itemize}
\item{}{\small
we have only to keep in mind
that to find a~tangent means to draw the line
that connects two points of the curve at an infinitely
small distance, or the continued side of a~polygon
with an~infinite number of angles which for us takes the place of the
curve.}
\end{itemize}

\noindent
We may conclude that Leibniz rested his calculus on
appealing to the structure of a curve ``in the small.''

\medskip
\begin{center}
\includegraphics{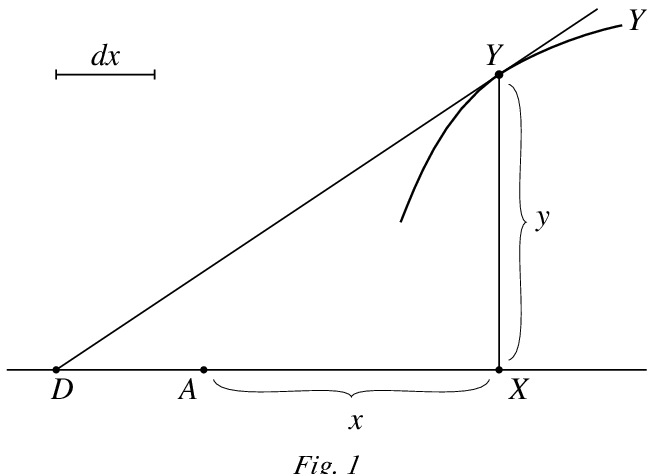}
\end{center}
\medskip

At that time, there were practically two standpoints as
regards the status of infinitesimals.
According to one of them, which seemed to be shared by
Leibniz,
an~infinitely small quantity was thought
of as an entity ``smaller than any given or assignable magnitude.''
Actual ``indivisible''  elements
comprising numerical quantities and geometrical figures are
the perceptions
corresponding to this concept of the infinitely small.
Leibniz did not doubt the existence of
``simple substances incorporated
into the structure of complex substances,''  i.e., {\it monads}.
``It is these monads that are the genuine atoms of nature
or, to put it short, elements of things'' \cite[p.~413]{Leib1}.

For the other founder of analysis, Newton,
the concept of infinite smallness is primarily related to the idea of
vanishing quantities
\cite{Newton, Reader}.
He viewed the indeterminate quantities
``not as made up of indivisible particles but as described by
a~ continuous motion'' and ``as increasing or decreasing by a~
perpetual motion, in their nascent or evanescent state.''

The celebrated
``method of prime and ultimate ratios''
reads in his classical treatise  {\it Mathematical Principles of
Natural Philosophy} (1687) as follows (see \cite[p.~101]{Reader}:

\begin{itemize}
\item{}{\small
Quantities, and the ratios of quantities, which in any finite
time converge continuously to equality, and before the end
of that time approach nearer to each other than by any given
difference, become ultimately equal.}
\end{itemize}

\noindent
Propounding the ideas which are nowadays attributed to
 the theory of limits, Newton
exhibited the insight, prudence, caution, and wisdom
characteristic of any great scientist  pondering over
the concurrent views and opinions.
He wrote (see \cite[p.~169]{Newton}):

\begin{itemize}
\item{}{\small
To institute an~analysis after
this manner in finite quantities and investigate the prime
or ultimate ratios of these finite quantities when
in their nascent state is consonant to the geometry
of the ancients, and I was willing to show that
in the method of fluxions there is no necessity of
introducing  infinitely small figures into geometry.

Yet the analysis may be performed in any kind of
figure,
whether finite or infinitely small, which are imagined
similar to the evanescent figures, as likewise in the
figures, which, by the method of indivisibles,
used
to be reckoned as infinitely small provided you proceed with
due caution.}
\end{itemize}

\noindent
Leibniz's
views were as much pliable and in-depth
dialectic. In his  famous  letter to Varignion
of February 2, 1702 \cite{Reader}, stressing the idea that
``it is unnecessary to make
mathematical analysis
depend on metaphysical controversies,''
he pointed out  the unity of the concurrent views
of the objects of the new calculus:

\begin{itemize}
\item{}{\small
If any opponent tries to contradict this proposition,
it follows from our calculus that the error will be less
than any possible assignable error, since it is in our
power to make this incomparably small magnitude small
enough for this purpose, inasmuch as we can always take
a~magnitude as small as we  wish.
Perhaps this is what you mean, Sir, when you speak on the
inexhaustible, and the rigorous demonstration of the
infinitesimal calculus which we use undoubtedly is to
be found here.\,...

So it can also be said that infinites
and infinitesimals are grounded in such a~way that  everything
in geometry, and even in nature, takes place as if they
were perfect realities. Witness not only our geometrical
analysis of transcendental curves but also my law
of continuity, in virtue of which it is permitted to consider
rest as infinitely small motion (that is,
as equivalent
to a~
species of its own contradictory),
and coincidence as infinitely small distance, equality
as the last inequality, etc.}
\end{itemize}

\noindent
Similar views were expressed by Leibniz
in the following
quotation (see \cite{Leib1}[ p.~190]) whose end in italics
is often cited in works on infinitesimal analysis
in the wake of Robinson
\cite[pp.~260--261]{ROB}:

\begin{itemize}
\item{}{\small
There is no need to take the infinite here rigorously,
but only as when we say in optics that the rays
of the sun come from a point infinitely distant, and
thus are regarded as parallel.
And when there are more degrees of infinity,
or infinitely small, it is as the sphere of the earth
is regarded as a point in respect
to the distance of the sphere of the fixed stars,
and a ball which we hold in the hand is also a point
in comparison with the semidiameter of the sphere of the earth.
And then the distance to the fixed stars is infinitely infinite
or an infinity of infinities in relation to the diameter of the ball.
For in place of the infinite or the infinitely small we
can take quantities as great or as small as is necessary in order
that the error will be less than any given error.
In this way {\it we only differ from the style of Archimedes
in the expressions, which are more direct in our method and
better adapted to the art of discovery}}.
\end{itemize}

\section{L.~Euler}

The 18th century is rightfully
called the age of Euler
in the history of mathematical analysis (cf.~\cite{Boyer}).
Everyone looking through his textbooks \cite{Omni}
will be staggered by  subtle technique and in-depth
penetration into the essence of the subject.
It is worth recalling that an~outstanding Russian engineer and scientist
Krylov went into raptures at the famous
Euler formula
$e^{i\pi}=-1$
viewing it as
the~quintessential symbol of integrity of all branches of mathematics.
He noted in particular that ``here $1$
presents arithmetic;  $i$, algebra;
$\pi$, geometry; and  $e$, analysis.''

Euler demonstrated an~open-minded approach,
which  might deserve the epithet ``systemic''  today,
to studying
mathematical problems: he applied the most sophisticated tools
of his time.
We must observe that part and parcel of his research was
the effective and productive use of various infinitesimal
concepts, first of all, infinitely large and
infinitely small numbers.
Euler
thoroughly explained  the
methodological background of his technique
in the form of the ``calculus of zeros.''
It is a popular fixation  to
claim that nothing is perfect  and to enjoy the imaginary failures
and
follies of the men of genius (``to look for sun-spots'' in the words of a~Russian saying).
For many years Euler had been incriminated
in the ``incorrect''  treatment of divergent series until
his ideas were fully accepted at the turn of the 20th century.
We may encounter such a~phrase in the literature:
``As to the problem of divergent series, Euler
was sharing quite an~up-to-date point of view.''
It would be more fair to topsy-turvy this phrase and say
that the mathematicians of today have finally caught  up with
some of Euler's ideas.
In fact the opinion that ``we cannot admire the way Euler
corroborates his analysis by introducing zeros of various orders''
is as self-assured as the statement that
``the giants of science, mainly  Euler and Lagrange, have
laid false foundations of analysis.''
It stands to reason to admit once and for ever that Euler was in full
possession of analysis and completely aware what he had created.

\section{G.~Berkeley}

The general ideas of analysis greatly
affected the lineaments of the ideological outlook in the 18th century.
The most vivid examples of the depth of penetration of the notions of
infinitely large and infinitely small quantities into the
cultural media of that time are  in particular
{\it Gulliver's Travels\/} by Jonathan Swift published in 1726
(Lilliput and Brobdingnag) and the celebrated {\it Micromegas 1752\/}
written by bright and venomous F.~M.~Arouer, i.e., Voltaire.
Of interest is the fact that as an~epigraph for his classical
treatise \cite{ROB}, Robinson
chose  the beginning of the following
speech of Micromegas (cf.~\cite[p.~154]{Volt}):

\begin{itemize}
\item{}{\small
Now I see clearer than ever that nothing can be judged by its
visible magnitude.
Oh, my God, who granted reason to creatures of such tiny sizes!
An~infinitely small thing is equal to an~infinitely large
one when facing you; if living beings still smaller than those
were possible, they could have reason exceeding the
intellect of those magnificent creatures of yours which
I can see in the sky, and one foot of which could cover the
earth.}
\end{itemize}

\noindent
 A~serious  and dramatic impact on the development of
infinitesimal  analysis was made in 1734 by
Bishop Berkeley,
a~great cleric and theologian, who published the~pamphlet
{\it The Analyst, or a~Discourse Addressed to an~
Infidel Mathematician, wherein it is examined whether the object,
principles and inferences of the modern analysis are more
deduced than religious mysteries and points of faith\/}
\cite{Berkl}. By the way, this Infidel Mathematician
was E.~Halley, a brilliant astronomer and
a~ young friend of Newton.
The clerical spirit of this article by Berkeley
is combined with aphoristic observations and
killing precision of expression. The leitmotif of his
criticism of analysis reads:
``Error may bring forth truth, though it cannot
bring forth science.''

\noindent
Berkeley's
challenge was addressed to all natural sciences:

\begin{itemize}
\item{}{\small
I have no controversy about your conclusions, but only
about your logic and method.
How do you demonstrate? What objects are you conversant with, and
whether you conceive them clearly?
What principles you proceed upon; how sound they may be; and
how you apply them?}
\end{itemize}

\noindent
Berkeley's invectives could not  be left unanswered by
the most progressive representatives of the scientific thought
of the 18th century, the encyclopedists.

\section{J.~D'Alembert and L.~Carnot}

A~turning point in the history of  the basic notions
of analysis is associated with the ideas and activities
of D'Alembert,
one of the initiators and leading authors
of the immortal masterpiece of the thought of the
Age of Enlightenment, the French
{\it Encyclopedia or Explanatory  Dictionary of Sciences, Arts, and Crafts}.

In the article ``Differential'' he wrote:
``Newton never considered differential calculus to be
some calculus of the infinitely small, but he rather viewed
it as the method of prime and ultimate ratios''
\cite[p.~157]{Reader}.
D'Alembert was the first mathematician who declared that
he had found the proof that the infinitely small
``do exist neither in Nature nor in the assumptions
of geometricians''
(a~quotation from his article ``Infinitesimal'' of 1759).

 The D'Alembert standpoint in
{\it Encyclopedia\/}  contributed much to the
formulation by the end of the 18th century of the
understanding of an infinitesimal as a~
vanishing magnitude.

It seems worthy to recall in this respect
the~book by Carnot
{\it Considerations on Metaphysics
of the Infinitely Small\/}
wherein he observed that
``the notion of infinitesimal is less clear than that
of limit implying nothing else but the
difference between such a~limit and the quantity
whose
ultimate value it provides.''

\section{B.~Bolzano, A.~Cauchy, and  K.~Weierstrass}

The 19th century was the time of  building analysis
over the theory of limits.
Outstanding contribution to this process belongs
to Bolzano, Cauchy, and
Weierstrass
whose achievements are mirrored in every traditional textbook
on differential and integral calculus.

The new canon of rigor  by Bolzano, the definition
by Cauchy of an~infinitely small quantity
as a~vanishing variable and, finally,
the $\varepsilon$-$\delta$-technique
by
Weierstrass are indispensable to the history
of mathematical thought, becoming part and parcel of the modern culture.

It is worth observing   (see \cite{Reader}) that,
giving a~verbal definition of continuity,
both Cauchy and Weierstrass chose practically
the same words:

\medskip
\hskip.5cm\vbox
{\hangindent = 20 pt
{\small An~infinitely small increment given to the~variable \hfill\break
\phantom{}produces an infinitely small increment of the function itself.
\hfill\break
\phantom{MMMMMMMMMMMMMMMMMMMMMMMMMM}{\eightsc Cauchy}}\hfill
}

\medskip
\hskip.5cm\vbox
{\hangindent = 20 pt
{\small Infinitely small variations in the arguments\hfill\break
\phantom{}correspond to those of the function.
\phantom{Jennnnnnnnnnnnnnnnnnnnn}\hfill\break
\phantom{MMMMMMMMMMMMMMMm}{\small\eightsc Weierstrass}}\hfill
}
\medskip

\noindent
This coincidence  witnesses the respectful desire
of the  noble authors  to interrelate the
new ideas with the views of their great predecessors.

Speculating about significance of the change of
analytical views in the 19th century, we should always
bear in mind the important observation by Severi
\cite[p.~113]{Sev} who wrote:

\begin{itemize}
\item{}{\small
This reconsideration, close to
completion nowadays, has however no ultimate value most
scientists believe in.
Rigor itself is, in fact, a~function of the amount
of knowledge at each historical period, a function that
corresponds to the manner in which science handles the truth.}
\end{itemize}

\section{N.~N.~Luzin}

The beginning of the 20th century in mathematics
was marked by a growing distrust of the concept of infinitesimal.
This tendency became prevailing as mathematics was
reconstructed on the set-theoretic  foundation
whose proselytes gained the key strongholds in the  1930s.

In the first edition of the {\it Great Soviet Encyclopedia\/}
in 1934, Luzin wrote (cf.~ \cite[pp.~293--294]{Luz59}):

\begin{itemize}
\item{}{\small
As to a~constant infinitely small quantity other than zero,
the modern mathematical analysis, without discarding the formal
possibility of defining the idea of a~constant
infinitesimal (for instance, as a~corresponding segment in
some non-Archimedean geometry), views this idea as absolutely fruitless
since it turns out impossible to introduce such
an~infinitesimal into calculus.}
\end{itemize}

\noindent
The publication  of the~textbook
{\it Fundamentals of Infinitesimal Calculus\/}
by
Vygodski\u\i{}
became a~noticeable event  in Russia at that time and
gave rise to a~serious and
sharp criticism. Vygodski\u\i{} tried to preserve the concept
of infinitesimal by appealing to  history and paedeutics.

\noindent
He wrote in particular (cf.~\cite[p.~160]{Vyg}):

\begin{itemize}
\item{}{\small
If it were only the problem of creating some logical apparatus
that could work by itself then, having eliminated infinitesimals
from considerations and having driven differentials out
of mathematics, one could celebrate a~victory over the
difficulties that have been impeded the way of
mathematicians and philosophers during the last two centuries.
Infinitesimal  analysis originated  however from practical needs,
its relations with the natural sciences and technology
(and, later, with social sciences) becoming increasingly
strong and fruitful in the course of time.
Complete elimination of infinitesimals would
hinder these relations or even make them impossible.}
\end{itemize}

\noindent
Discussing this textbook by Vygodski\u\i,
Luzin wrote in the 1940s (cf.~\cite[p.~398]{Luz59}):

\begin{itemize}
\item{}{\small
This course, marked by internal  integrity and lit by
the great idea the author remains faithful to,
falls beyond the framework of the style in which the modern
mathematical analysis has been developed for 150 years
and which is now nearing its completion.}
\end{itemize}

\noindent
Luzin's attitude to infinitesimals deserves
special attention as apparent manifestation
and convincing evidence of the~background drama typical of the history of every
profound idea that
enchants and inspires the mankind.
Luzin
had a~unique capability of penetration into the essence
of the most intricate mathematical problems, and
he might be said to possess a~remarkable gift of foresight
\cite{Lavr79, Lavr80, Luz84}.

The concept of actual infinitesimals seemed to be extremely
appealing to him psychologically, as he emphasized \cite[p.~398]{Luz59}:

\begin{itemize}
\item{}{\small
The idea about them has never been successfully
driven out of my mind.
There are obviously some deeply hidden reasons still unrevealed
completely that make our minds inclined to looking
at
infinitesimals favorably.}
\end{itemize}

\noindent
In  one of his letters to Vygodski\u\i\ which was written in 1934
he predicted that ``infinitesimals will be
fully rehabilitated from
a perfectly  scientific point of view as kind of `mathematical
quanta.'\,''
In another of his publications (cf.~\cite{Vil}), Luzin sorrowfully remarked:

\begin{itemize}
\item{}{\small
When the mind starts  acquaintance with analysis, i.e.,
during the mind's spring season, it is the infinitesimals,
which deserve to be called the ``elements'' of quantity,
that the mind begins with.
However, surfeiting itself
gradually
with knowledge, theory, abstraction and fatigue,
the mind  gradually forgets its primary intentions,
smiling at their ``childishness.''
In short, when the mind is in its autumn season,
it allows itself to become convinced of the unique
sound foundation by means of limits.}
\end{itemize}

\noindent
This limit conviction was energetically corroborated
by Luzin
in his textbook
{\it Differential Calculus\/}   wherein he  particularly emphasized
\cite[p.~61]{Luz61}:

\begin{itemize}
\item{}{\small
To  grasp the very
{\it essence of the matter\/}
correctly,
the student should first of all made it clear that each infinitesimal
is always
a~variable quantity
{\it by its very definition};
therefore, no constant number,
however tiny, is {\it ever\/} infinitely small.
The student should beware of using comparisons or
similes of such a kind for instance as
``One centimeter is a~magnitude infinitely small as compared with
the diameter of the sun.''  This phrase is pretty incorrect.
Both magnitudes, i.e., one centimeter and the diameter
of the sun, are constant quantities and so they are
{\it finite},
one  much smaller than the other, though.
Incidentally, one centimeter is not a~small length at all as
compared for instance with the ``thickness of a~hair,''
becoming  a long distance  for a~moving microbe.
In order to eliminate any risky comparisons and
haphazard subjective similes, the student
{\it must remember that neither constant magnitude is infinitesimal
nor any number, however small  these might be}.
Therefore, it would be quite appropriate to abandon the term
``{\it infinitesimal magnitude\/}''  in favor of  the term
``{\it infinitely vanishing variable},''
as the latter expresses the {\it idea of variability\/}
most vividly.}
\end{itemize}

\section{A.~Robinson}

The seventh posthumous edition of this textbook by
Luzin
was published in 1961
simultaneously with  Robinson's
{\it Non-Standard Analysis\/}  which laid a~modern foundation
for the calculus of infinitesimals.
Robinson based his research on the local theorem by
Mal$'$tsev,
stressing its
``fundamental importance for our theory'' \cite[p.~13]{ROB}
and giving explicit references  to
Mal$'$tsev's
article
dated as far back as 1936.
Robinson's discovery  elucidates the ideas of the
founders of differential and integral calculus,
witnessing  the spiral evolution of mathematics.

\bibliographystyle{plain}

\enddocument